\documentclass[a4paper,12pt]{article}

\usepackage{amssymb}
\usepackage{graphicx}
\usepackage{amsmath}
\usepackage{amsfonts}
\usepackage{thmtools}

\newcommand{\singlespacing}{\let\CS=\@currsize\renewcommand{\baselinestreatch}{1.0}\tiny\CS}
\newcommand{\doublespacing}{\let\CS=\@currsize\renewcommand{\baselinestreatch}{1.5}\tiny\CS }

\newtheorem{theorem}{Theorem}[section]
\newtheorem{corollary}{Corollary}[section]

\newtheorem{lemma}{Lemma}[section]

\newtheorem{remark}{Remark}[section]

\numberwithin{equation}{section}

\begin{document}
\begin{center}
{\Large \bf{The $k$-almost Ricci solitons and contact geometry}}\\
\vspace{.3cm}
Amalendu Ghosh $^1$ and Dhriti Sundar Patra $^2$
\end{center}
\newtheorem{Theorem}{\quad Theorem}[section]
\newtheorem{Definition}[Theorem]{\quad Definition}
\newtheorem{Corollary}[Theorem]{\quad Corollary}
\newtheorem{Lemma}[Theorem]{\quad Lemma}
\newtheorem{Example}[Theorem]{\emph{Example}}
\newtheorem{Proposition}[Theorem]{Proposition}
\numberwithin{equation}{section}
{\textbf{Abstract:} The aim of this article is to study the $k$-almost Ricci soliton and $k$-almost gradient Ricci soliton on contact metric manifold. First, we prove that  if a compact $K$-contact metric is a $k$-almost gradient Ricci soliton then it is isometric to a unit sphere $S^{2n+1}$. Next, we extend this result on a compact $k$-almost Ricci soliton when the flow vector field $X$ is contact. Finally, we study some special types of $k$-almost Ricci solitons where the potential vector field $X$ is point wise collinear with the Reeb vector field $\xi$ of the contact metric structure.}

\noindent\\
\textbf{Mathematics Subject Classification 2010}: 53C25, 53C20, 53C15\\\\
\textbf{Keywords}: Contact metric manifold, $k$-almost Ricci soliton, $k$-almost gradient Ricci soliton, $K$-contact manifold, Sasakian manifold, Einstein manifold.

\section{Introduction}
A Riemannian manifold $(M^n,g)$ is said to be a  Ricci soliton if there exists a vector field $X$ on $M^n$ and a constant $\lambda$ satisfying the equation $S + \frac{1}{2}\pounds_{X}g = \lambda g$, where $\pounds_{X}g$ denotes the Lie-derivative of $g$ along the vector field $X$ on $M^n$ and $S$ is the Ricci tensor of $g$. In general, $X$ and $\lambda$ are known as the potential vector field and the soliton constant, respectively. Ricci solitons are the fixed points of Hamilton's Ricci flow: $\frac{\partial}{\partial t}g(t) = -2 S(g(t))$ (where $g(t)$ a one-parameter family of metrics on $M^n$) viewed as a dynamical system on the space of Riemannian metrics modulo diffeomorphisms and scalings (cf. \cite{HRS}). Recently, the notion of Ricci soliton was generalized by Pigoli-Rigoli-Rimoldi-Setti \cite{PRRS} to almost Ricci soliton by allowing the soliton constant $\lambda$ to be a smooth function.\\

Recently, Wang-Gomes-Xia \cite{WGX} extended the notion of almost Ricci soliton to $k$-almost Ricci soliton which is defined as

\begin{Definition}
A complete Riemannian manifold $(M^n,g)$ is said to be a $k$-almost Ricci soliton, denoted by $(M^n, g, X, k, \lambda)$, if there exists a smooth vector field $X$ on $M^n$, a soliton function $\lambda\in C^\infty(M^n)$ and a positive real valued function $k$ on $M^n$ such that
 \begin{eqnarray}\label{1.1}
&S + \frac{k}{2}\pounds_{X}g = \lambda g.
 \end{eqnarray}
\end{Definition}

This notion has been justified as follows. Suppose $(M^n,g_{0})$ be a complete Riemannian manifold of dimension $n$ and let $g(t)$ be a solution of the Ricci flow equation defined on $[0,\epsilon)$, $\epsilon > 0$, such that $\psi_{t}$ is a one-parameter family of diffeomorphisms of $M^n$, with $\psi_{0} = id_{M}$ and $g(t)(x)  = \rho(x,t)\psi_{t}^* g_{0}(x)$ for every $x\in M^n$, where $\rho(x,t)$ is a positive smooth function on $M^n\times [0,\epsilon)$. Then one can deduce
\begin{eqnarray*}
&\frac{\partial}{\partial t}g(t)(x) = \frac{\partial}{\partial t}\rho(x,t) \psi_{t}^* g_{0}(x) + \rho(x,t)\psi_{t}^* \pounds_{\frac{\partial}{\partial t}\psi(x,t)}g_{0}(x).
\end{eqnarray*}
When $t = 0$, the foregoing equation reduces to
\begin{eqnarray*}
&S_{g_{0}} + \frac{k}{2}\pounds_{X}g_{0} = \lambda g_{0},
\end{eqnarray*}
where $X = \frac{\partial}{\partial t}\psi(x,0)$, $\lambda(x) = -\frac{1}{2} \frac{\partial}{\partial t}\rho(x,0)$ and $k(x) = \rho(x,0)$.\\

A $k$-almost Ricci soliton is said to be $shrinking$, $steady$ or $expanding$ accordingly as $\lambda$ is positive, zero or negative, respectively. It is trivial (Einstein) if the flow vector field $X$ is homothetic, i.e. $\pounds_{X}g = cg$, for some constant $c$. Otherwise, it is non-trivial. A $k$-almost Ricci soliton is said to be a $k$-almost gradient Ricci soliton if the potential vector field $X$ can be expressed as a gradient of a smooth function $u$ on $M^n$, i.e. $X = Du$, where $D$ is the gradient operator of $g$ on $M^n$. In this case, we denotes $(M^n,g,Du,k,\lambda)$ as a $k$-almost gradient Ricci soliton with potential function $u$. Further, the fundamental equation (\ref{1.1}) takes the form
\begin{eqnarray}\label{1.2}
&S + k\nabla^2 u = \lambda g,
\end{eqnarray}
where $\nabla^2 u$ denotes the Hessian of $u$.\\

In particular, a Ricci soliton is the $1$-almost Ricci soliton with constant $\lambda$, and an almost Ricci soliton is just the $1$-almost Ricci soliton. Barros and Ribeiro Jr. proved (cf. \cite{BR}) that a compact almost Ricci soliton with constant scalar curvature is isometric to a Euclidean sphere. An analogous result has also been proved by Wang-Gomes-Xia \cite{WGX} for the case of $k$-almost Ricci soliton.\\

\textbf{Theorem [WGX]:} Let $(M^n, g, X, k, \lambda)$, $n \geq 3$ be a non-trivial $k$-almost Ricci soliton with constant scalar curvature $r$. If $M^n$ is compact, then it is isometric to a standard sphere $S^n(c)$ of radius $c = \sqrt{\frac{2n(2n + 1)}{r}} $.\\

Recall that a smooth manifold $M^n$ together with a Riemannian metric $g$ is said to be a generalized quasi-Einstein manifold if there exist smooth functions $f$, $\mu$ and $\lambda$ such that (cf. \cite{Cat})
\begin{align*}
S + \nabla^2f - \mu df\otimes df = \lambda g.
\end{align*}
For $\mu = \frac{1}{m}$, the generalized quasi-Einstein manifold is known as generalized $m$-quasi-Einstein manifold (cf. \cite{Bar-Gom, Bar-Rib}), and when $\lambda$ is constant the generalized quasi-Einstein manifold is simply known as $m$-quasi-Einstein manifold. Case-Shu-Wei \cite{CSW} proved that any complete quasi-Einstein-metric with constant scalar curvature is trivial (Einstein). Subsequently, this has been extended by Barros-Gomes \cite{Bar-Gom}. In fact, they proved that any compact generalized $m$-quasi-Einstein metric with constant scalar curvature is isometric to a standard Euclidean sphere $S^n$. Particularly, Barros-Ribeiro \cite{Bar-Rib} construct a family of nontrivial generalized $m$-quasi-Einstein metric on the unit sphere $S^n(1)$ that are rigid in the class of constant scalar curvature. It is interesting to note that by a suitable choice of the function $f$ it is possible to reduce any generalized $m$-quasi-Einstein metric to a $k$-almost Ricci soliton. For instance, if we take $u=e^{\frac{f}{m} }$ and $k=-\frac{m}{u}$, then (\ref{1.2}) reduces to
 \begin{align*}
S + \nabla^2f - \frac{1}{m} df\otimes df = \lambda g.
\end{align*}
Thus, in one hand $k$-almost Ricci soliton generalizes generalized $m$-quasi-Einstein metric, on the other it covers gradient Ricci soliton and gradient almost Ricci soliton.  For details we refer to \cite{WGX}. Recently, Yun-Co-Hwang \cite{YCH} studied Bach-flat $k$-almost gradient Ricci solitons.\\

During the last few years Ricci soliton and almost Ricci soliton have been studied by several authors (cf. \cite{S}, \cite{GSC}, \cite{CS}, \cite{GS} and \cite{AGC}) within the frame-work of contact geometry. In \cite{S}, Sharma initiated the study of gradient Ricci soliton within the frame-work of $K$-contact manifold and prove that \textit{``any complete $K$-contact metric admitting a gradient Ricci soliton is Einstein and Sasakian"}. Later on, this has been generalized by the second author \cite{AGC} who proved that \textit{``if a complete $ K $-contact metric (in particular, Sasakian) represents a gradient almost Ricci soliton, then is it isometric to the unit sphere $S^{2n+1}$"}. Inspired by these results, here we consider contact metric manifolds whose metric is a $k$-almost Ricci soliton. Following \cite{Bar-Rib}, one can construct a family of non-trivial examples of generalized $m$-quasi-Einstein metrics on the odd dimensional unit sphere $S^{2n+1}$. Another motivation arises from the fact that any odd dimensional unit sphere satisfies the generalized $m$-quasi-Einstein condition and hence it satisfies the gradient $k$-almost Ricci soliton equation (\ref{1.2}). Since any odd dimensional unit sphere $S^{2n+1}$ admits standard $ K $-contact (Sasakian) structure we are interested in studying $K$-contact metric as a gradient $k$-almost Ricci soliton. We address this issue in Section 3 and prove that if a compact $K$-contact manifold admits a $k$-almost gradient Ricci soliton then it is isometric to a unit sphere $S^{2n + 1}$. Next, we study $k$-almost Ricci soliton in the frame-work of compact $K$-contact manifold when the potential vector field is contact. Finally, a couple of results on contact metric manifolds admitting $k$-almost Ricci soliton are presented under the assumption that the potential vector field $X$ is point wise collinear with the Reeb vector field $\xi$ of the contact metric structure.

\section{Preliminaries}First, we recall some basic definitions and formulas on a contact metric manifold. By a contact manifold we mean a Riemannian manifold $M^{2n+1}$ of dimension $(2n + 1)$ which carries a global $1$-form $\eta$ such that $\eta\wedge(d\eta)^{n}\ne0$  everywhere on $M^{2n+1}$. The form $\eta$ is usually known as the contact form on $M^{2n+1}$. It is well known that a contact manifold admits an almost contact metric structure on $(\varphi,\xi,\eta, g)$, where $\varphi$ is a tensor field of type $(1,1)$, $\xi$ a global vector field known as the characteristic vector field (or the Reeb vector field) and $g$ is Riemannian metric, such that
\begin{eqnarray}
&\varphi^2 Y = -Y + \eta(Y) \xi,  \label{2.1C} \\
&\eta(Y) = g(Y, \xi),  \label{2.1B} \\
&g(\varphi Y, \varphi Z) = g(Y,Z) - \eta(Y)\eta(Z), \label{2.2}
\end{eqnarray}
for all vector fields $Y$, $Z$ on $M$. It follows from the above equations that $\varphi\xi=0$ and $\eta \circ\varphi=0$ (see \cite{Blair}, p.43). A Riemannian manifold $M^{2n+1}$ together with the almost contact metric structure $(\varphi,\xi,\eta,g)$ is said to be a contact metric if it satisfies ( \cite{Blair}, p.47)
\begin{eqnarray}\label{2.1A}
&d\eta(Y,Z)=g(Y,\varphi Z),
\end{eqnarray}
for all vector fields $Y$, $Z$ on $M$. In this case, we say that $g$ is an associated metric of the contact metric structure. On a contact metric manifold $M^{2n+1}(\varphi,\xi,\eta,g),$ we consider two self-adjoint operators $h=\frac{1}{2}\pounds_{\xi}\varphi$ and $l=R(.,\xi)\xi$, where $\pounds_{\xi}$ is the Lie-derivative along $\xi$ and $R$ is the Riemann curvature tensor of $ g $. The two operators $h$ and $l$ satisfy (e.g., see \cite{Blair}, p.84, p.85)
$$Tr~h=0, \hskip 0.2cm Tr~(h\varphi)=0, \hskip 0.2cm h\xi=0, \hskip 0.2cm l\xi = 0, \hskip 0.2cm h\varphi=-\varphi h.$$
We now recall the following
\begin{lemma}(\cite{Blair}, p.84; p.112; p.111)
On a contact metric manifold $M^{2n+1}(\varphi, \xi, \eta, g)$ we have
\begin{eqnarray}
&\nabla_{Y}\xi=-\varphi Y-\varphi hY, \label{2.3}\\
&Ric_{g}(\xi,\xi)=g(Q\xi,\xi)=Tr~l=2n-Tr~(h^{2}), \label{2.4}\\
&(\nabla_{Z}\varphi)Y + (\nabla_{\varphi Z}\varphi)\varphi Y = 2g(Y,Z)\xi-\eta(Y)(Z + hZ + \eta(Z)\xi), \label{2.5}
\end{eqnarray}
for all vector fields $Y$, $Z$ on $M$; where $\nabla$ is the operator of covariant differentiation of $g$ and $Q$ the Ricci operator associated with the $(0,2)$ Ricci tensor given by $S(Y,Z)=g(QY,Z)$ for all vector fields $Y$, $Z$ on $M$.
\end{lemma}
A contact metric manifold is said to be $K$-contact if the vector field $\xi$ is Killing, equivalently if $h=0$ (\cite{Blair}, p.87). Hence on a $K$-contact manifold Eq$.$ (\ref{2.3}) becomes
\begin{eqnarray}\label{2.6}
\nabla_{Y}\xi = - \varphi Y,
\end{eqnarray}
for any vector field $Y$ on $M$. Moreover, on a $K$-contact manifold the following formulas are also valid.
\begin{lemma}(see Blair \cite{Blair}, p.113; p.116)
On a $K$-contact manifold \\$M^{2n+1}(\varphi, \xi, \eta, g)$ we have
\begin{eqnarray}
&Q \xi = 2n \xi,  \label{2.7} \\
&R(\xi,Y)Z=(\nabla_{Y}\varphi)Z, \label{2.9}
\end{eqnarray}
for all vector fields $Y$, $Z$ on $M$.
\end{lemma}
An almost contact metric structure on $M$ is said to be normal if
the almost complex structure $J$ on $M\times \mathbb{R}$ defined by (e.g., see Blair \cite{Blair}, p.80)
\begin{eqnarray*}
&J(X,fd/dt)=(\varphi X-f\xi,\eta(X)d/dt),
\end{eqnarray*}
where $f$ is a real function on $M\times R$, is integrable. A normal contact metric manifold is said to be Sasakian. On a Sasakian manifold (e.g., \cite{Blair}, p.86)
\begin{eqnarray*}
(\nabla_{X}\varphi)Y=g(X,Y)\xi-\eta(Y)X,
\end{eqnarray*}
for all vector fields $X$, $Y$ on $M$. Further, a contact metric manifold is Sasakian if and only if the curvature tensor $R$ satisfies (e.g., \cite{Blair}, p.114)
\begin{eqnarray}\label{2.10}
R(X,Y)\xi=\eta(Y)X-\eta(X)Y,
\end{eqnarray}
for all vector fields $X$, $Y$ on $M$. A Sasakian manifold is $K$-contact but the converse is true only in dimension $3$ (e.g., \cite{Blair}, p.87).\\

A contact metric manifold is said to be  $\eta$-Einstein if the Ricci tensor $S$ satisfies $S(Y,Z) = a g(Y,Z) + b\eta(Y)\eta(Z)$
for any vector fields $Y$, $Z$ on $M$ and are arbitrary functions $a$, $b$ on $M$. The functions $a$ and $b$ are constant for a $K$-contact manifold of dimension $> 3$ (cf. \cite{YK}).\\

Let $T^{1}M$ be the unit tangent bundle of a compact orientable Riemannian manifold $(M,g)$ equipped with the Sasaki metric $g_{s}$. Any unit vector field $U$ determines a smooth map between $(M,g)$ and $(T^{1}M, g_{s})$. The energy $E(U)$ of the unit vector field $U$ is defined by
\begin{eqnarray*}
&E(U) = \frac{1}{2}\int\parallel dU \parallel^{2} dM = \frac{n}{2}vol(M, g) + \frac{1}{2}\int_{M}\parallel \nabla U \parallel^{2} dM,
\end{eqnarray*}
where $dU$ denotes the differential of the map $U$ and $dM$ denotes the volume element of $M$. $U$ is said to be a harmonic vector field if it is a critical point of the energy functional $E$ defined on the space $\chi^{1}$ of all unit vector fields on $(M,g)$. A contact metric manifold is said to be an $H$-contact manifold if the Reeb vector field $\xi$ is harmonic. In \cite{Perrone}, Perrone proved that \textit{``A contact metric manifold is an $H$-contact manifold, that is $\xi$ is a harmonic vector field, if and only if $\xi$ is an eigenvector of the Ricci operator."}
On a contact metric manifold, $\xi$ is an eigenvector of the Ricci operator implies that $ Q\xi = (Tr~l)\xi $. This is true for many contact metric manifolds. Such as, $\eta$-Einstein contact metric manifolds, $ K $-contact (in particular Sasakian) manifolds, $(k,\mu)$-contact manifolds and the tangent sphere bundle of a Riemannian manifold of constant curvature. In particular, this condition holds on the unit sphere $S^{2n+1}$ with standard contact metric structure.

\begin{Definition}
A vector field $X$ on a contact manifold is said to be a contact vector field if it preserve the contact form $\eta$, i.e.
\begin{eqnarray}\label{2.11}
\pounds_{X}\eta = f\eta,
\end{eqnarray}
for some smooth function $f$ on $M$. When $f = 0$ on $M$, the vector field $X$ is called a strict contact vector field.
\end{Definition}

\textbf{Example 2.2: }
It is well know \cite{Blair} that any odd dimensional unit sphere $ S^{2n+1} $ admits a standard $ K $-contact (Sasakian) structure $(\varphi, \xi, \eta, g)$ and hence the Reeb vector field satisfies (\ref{2.6}), for any vector field $ Y $ on $ S^{2n+1} $. We now recall the theorem of Obata \cite{O} that {\it a complete connected Riemannian manifold $(M,g)$ of dimension $ > 2 $ is isometric to a sphere of radius $\frac{1}{c}$ if and only if it admits a non-trivial solution $ k $ of the equation $\nabla \nabla k = - c^{2}k g $}. For unit sphere this transforms to $\nabla \nabla k = - k g $, where $ k $ is the eigenfunction of the Laplacian on $ S^{2n+1} $. Let $ X $ be a vector field on $ S^{2n+1} $ such that $ X = -Dk + \mu \xi $, where $ \mu $ is a constant. Differentiating this along an arbitrary vector field $Y$ on $S^{2n+1}$ and using (\ref{2.6}) we obtain $\nabla_Y X = -\nabla_{Y}Dk - \mu\varphi Y$. Then by Obata's theorem and (\ref{2.6}) we see that
\begin{eqnarray}\label{2.14}
&\frac{k}{2}(\pounds_{X}g)(Y, Z) + S(Y, Z) = (k^2+2n)g(Y, Z),
\end{eqnarray}
for all vector fields $Y$, $Z$ on $ S^{2n+1} $. This shows that $(S^{2n+1}, g, X, \lambda)$ is a almost Ricci soliton with $ \lambda = k^2 + 2n $. Moreover, if we take $ X = Du $, for some smooth non constant function $u$ on $S^{2n+1}$, then from (\ref{2.14}) it follows that $S^{2n+1}$ also admits $k$-almost gradient Ricci soliton.

\section{ $K$-contact metric as $k$-almost gradient Ricci soliton and $k$-almost Ricci soliton }
We assume that a $K$-contact metric $g$ is a $k$-almost gradient Ricci soliton with the potential function $u$. Then the $k$-almost gradient Ricci soliton Eq$.$ (\ref{1.2}) can be exhibited as
\begin{eqnarray}\label{3.2}
k\nabla_{Y}Du + QY = \lambda Y,
\end{eqnarray}
for any vector field $Y$ on $M$; where  $D$ is the gradient operator of $g$ on $M$. Taking the covariant derivative of (\ref{3.2}) along an arbitrary vector field $Z$ on $M$ yields
\begin{eqnarray*}
k\nabla_{Z}\nabla_{Y}Du  =& \frac{1}{k}(Zk)(QY -\lambda Y)- (\nabla_{Z}Q)Y -Q(\nabla_{Z}Y) \nonumber\\
&+ (Z\lambda)Y + \lambda\nabla_{Z}Y,
\end{eqnarray*}
for any vector field $Y$ on $M$. Using this and (\ref{3.2}) in the well known expression of the curvature tensor $R(Y,Z) = [\nabla_{Y},\nabla_{Z}] - \nabla_{[Y,Z]},$ we can easily find out the curvature tensor which is given by
\begin{eqnarray}\label{3.1}
&k R(Y,Z)Du = \frac{1}{k}(Yk)(QZ -\lambda Z) - \frac{1}{k}(Zk)(QY -\lambda Y) \nonumber\\
&+ (\nabla_{Z}Q)Y - (\nabla_{Y}Q)Z + (Y\lambda)Z - (Z\lambda)Y,
\end{eqnarray}
for all vector fields $Y$, $Z$ on $M$.\\

Before entering into our main results we prove the following.
\begin{lemma}
On a $ K $-contact manifold $M^{2n+1}(\varphi, \xi, \eta, g)$, we have
\begin{eqnarray}\label{3.4}
\nabla_{\xi}Q = Q\varphi - \varphi Q.
\end{eqnarray}
\end{lemma}
\textbf{Proof:} Since $\xi$ is Killing on a $K$-contact manifold, we have $(\pounds_{\xi}Q)Y = 0$ for any vector field $Y$ on $M$. Taking into account (\ref{2.6}) it follows that
\begin{eqnarray*}
0 &=&\pounds_{\xi}(QY)-Q(\pounds_{\xi}Y)\\
&=& \nabla_{\xi}QY - \nabla_{QY}\xi-Q(\nabla_{\xi}Y)+Q(\nabla_{Y}\xi)\\
&=&(\nabla_{\xi}Q)Y + \varphi QY - Q \varphi Y,
\end{eqnarray*}
for any vector field $Y$ on $M$. This completes the proof. ~~~~~~~~~~~~~~~~~~$\square$

\begin{theorem}
Let $(M^{2n+1},g,Du,k,\lambda)$ be a $k$-almost gradient Ricci soliton with the potential function $u$. If $(M,g)$ is a compact $K$-contact manifold, then it is isometric to a unit sphere $S^{2n + 1}$.
\end{theorem}
\textbf{Proof: }Firstly, taking covariant differentiation of (\ref{2.7}) along an arbitrary vector field $Y$ on $M$ and using (\ref{2.6}), we get
\begin{equation}\label{3.3}
(\nabla_{Y}Q)\xi = Q\varphi Y-2n\varphi Y.
\end{equation}
Now, replacing $\xi$ instead of $Y$ in (\ref{3.1}) and making use of the $K$-contact condition \eqref{2.7}, (\ref{3.3}) and (\ref{3.4}), we get
\begin{eqnarray*}
&k R(\xi,Z)Du = (\frac{\lambda - 2n}{k})(Z  k)\xi + \frac{1}{k}(\xi k)(QZ - \lambda Z) - 2n\varphi Z  \nonumber\\
&+ \varphi QZ + (\xi\lambda)Z - (Z\lambda)\xi,
\end{eqnarray*}
for any vector field $Z$ on $M$. Scalar product of the last equation with an arbitrary vector field $Y$ on $M$ and using (\ref{2.9}), we obtain
\begin{eqnarray}\label{3.6}
&kg(( \nabla_{Z}\varphi) Y,Du) + (\frac{\lambda - 2n}{k})(Z k)\eta(Y) + (\xi \lambda - \frac{\lambda}{k}(\xi k))g(Y,Z)\nonumber\\
&+ \frac{1}{k}(\xi k)g(QY,Z) + 2n g(\varphi Y,Z) - g(Q\varphi Y,Z) - (Z\lambda)\eta(Y) = 0,
\end{eqnarray}
for any vector field $Z$ on $M$. Next, substituting $Y$ by $\varphi Y$ and $Z$ by $\varphi Z$ in \eqref{3.6} and using \eqref{2.1C}, $\eta o\varphi=0$ and $\varphi\xi=0$ provides
\begin{eqnarray*}
&kg((\nabla_{\varphi Z}\varphi) \varphi Y,Du)  + (\xi \lambda- \frac{\lambda}{k}(\xi k))\{g(Y,Z) - \eta(Y)\eta(Z)\} \nonumber\\
&+ \frac{1}{k}(\xi k)g(Q\varphi Y,\varphi Z) + 2n g(\varphi Y,Z) - g(\varphi QY,Z) = 0,
\end{eqnarray*}
for all vector fields $Y$, $Z$ on $M$. Adding the preceding Eq$.$ with (\ref{3.6}) and using (\ref{2.5}) (where $h=0$, as $M$ is $K$-contact) yields
\begin{eqnarray*}
&2\{k(\xi u) + (\xi\lambda) - \frac{\lambda}{k} (\xi k)\}g(Y,Z) + \{\frac{\lambda}{k}(\xi k) -(\xi\lambda) - k(\xi u)\}\eta(Y)\eta(Z) \nonumber\\
&+ \{(\frac{\lambda - 2n}{k})(Z k) - k(Zu) - (Z\lambda)\}\eta(Y) + \frac{1}{k}(\xi k) g(QY,Z)  \nonumber\\
&+ 4n g(\varphi Y,Z)  -  g(Q\varphi Y + \varphi QY,Z) + \frac{1}{k}(\xi k) g(\varphi QY,\varphi Z) = 0,
\end{eqnarray*}
for all vector fields $Y$, $Z$ on $M$. Anti-symmetrizing the foregoing equation provides
\begin{eqnarray*}
&\{(\frac{\lambda - 2n}{k})(Z k) - k(Zu) - (Z\lambda)\}\eta(Y) - 2 g(Q\varphi Y + \varphi QY,Z) \nonumber\\
&-\{(\frac{\lambda - 2n}{k})(Y k) - k(Yu) - (Y\lambda)\}\eta(Z) + 8n g(\varphi Y,Z) = 0,
\end{eqnarray*}
for all vector fields $Y$, $Z$ on $M$. Moreover, substituting  $Y$ by $\varphi Y$ and $Z$ by $\varphi Z$ in the last equation and using the $K$-contact condition \eqref{2.7}, \eqref{2.1C}, $\eta o\varphi=0$ and $\varphi\xi=0$ gives
\begin{eqnarray*}
&g(Q\varphi Y + \varphi QY,Z) = 4n g(\varphi Y,Z),
\end{eqnarray*}
for all vector fields $Y$, $Z$ on $M$. It follows from last Eq$.$ that
\begin{eqnarray}\label{3.8}
&Q\varphi Y + \varphi QY = 4n\varphi Y,
\end{eqnarray}
for any vector field $Y$ on $M$. Let $\{e_{i},\varphi e_{i},\xi \},i= 1,2,3,.....,n$, be an orthonormal $\varphi$-basis of $M$ such that $Qe_{i} = \sigma_{i}e_{i}$. Thus, we have $\varphi Qe_{i} = \sigma_{i}\varphi e_{i}.$
Substituting $e_{i}$ for $Y$ in (\ref{3.8}) and using the foregoing equation, we obtain $Q\varphi e_{i}  = (4n-\sigma_{i})\varphi e_{i}$. Using the $\varphi $-basis and (\ref{2.7}), the scalar curvature $r$ is given by
\begin{eqnarray*}
r &=& g(Q\xi,\xi) + \sum_{i=1}^{n}[g(Qe_{i},e_{i}) + g(Q\varphi e_{i},\varphi e_{i})]\\
&=& g(Q\xi,\xi) + \sum_{i=1}^{n}[\sigma_{i}g(e_{i},e_{i}) + (4n-\sigma_{i})g(\varphi e_{i},\varphi e_{i})]\\
&=& 2n(2n+1).
\end{eqnarray*}
Therefore, the scalar curvature $r$ is constant. As $M$ is compact, Theorem [WGX] shows that $M$ is isometric to $S^{2n+1}(c)$, where $c=\sqrt[]{\frac{2n(2n+1)}{r}}$ is the radius of the sphere. Since $r = 2n(2n + 1)$, we have $c = 1$. Hence, $M$ is isometric to a unit sphere $S^{2n + 1}$. This completes the proof.~~~~~~~~~~~~~~~~~~~~~~~~~$\square$
\begin{remark}
From the last theorem we see that any compact $K$-contact manifold $M$ admitting a gradient $k$-almost Ricci soliton is isometric to a unit sphere and hence of constant curvature $1$. Consequently, $M$ is Sasakian. Since $k$-almost Ricci soliton covers Einstein manifold, we may compare this as an extension of the odd dimensional Goldberg conjecture which states that any compact Einstein $K$-contact manifold is Sasakian. For details, we refer to Boyer-Galicki \cite{BG}.
\end{remark}
In particular, the above result is also true for complete Sasakian manifolds.
\begin{Corollary}
Let $(M^{2n+1},g,Du,k,\lambda)$ be a $k$-almost gradient Ricci soliton with the potential function $u$. If $(M,g)$ is a complete Sasakian manifold, then it is compact and isometric to a unit sphere $S^{2n + 1}$.
\end{Corollary}
\textbf{Proof: }On a Sasakian manifold the Ricci operator $Q$ and $\varphi$ commutes, i.e. $Q\varphi = \varphi Q$ (see \cite{Blair}, p.116). Using this in (\ref{3.8}) implies $Q\varphi Y = 2n \varphi Y$ for any vector field $Y$ on $M$. Substituting $Y$ by $\varphi Y$ in the last equation and using (\ref{2.7}) gives $QY = 2nY$ for any vector field $Y$ on $M$. This shows that $M$ is Einstein with Einstein constant $2n$. As $(M,g)$ is complete, $M$ is compact by Myers' Theorem \cite{MSB}. The rest of the proof follows from the last theorem.~~~~~~~~~~~~~~~~~~~~~~~~~~~~~~~~~~~~~~~~~~~~~~~~~~~~~~~~~~~~~~~~~~~~~~~~~~~~~~~~~~~~~~~~~~$\square$\\

Next, we extend the Theorem $3.1$  and consider $K$-contact metric as a $k$-almost Ricci soliton when its potential vector field is a contact vector field and prove
\begin{theorem}
Let $M^{2n+1}(\varphi,\xi,\eta,g)$ be a compact $K$-contact manifold with $X$ as a contact vector field. If $g$ is a $k$-almost Ricci soliton with  $X$ as the potential vector field, then $M$ is isometric to a unit sphere $S^{2n + 1}$.
\end{theorem}
\textbf{Proof: }
Firstly, taking Lie-derivative of (\ref{2.1A}) along $X$ and using (\ref{1.1}) we have
\begin{eqnarray}\label{4.5AB}
k(\pounds_{X}d\eta)(Y,Z) = 2 g(-QY+\lambda Y,\varphi Z) + kg(Y,(\pounds_{X}\varphi)Z),
\end{eqnarray}
for all vector fields $Y$, $Z$ on $M$. As $X$ is a contact vector field, we deduce from \eqref{2.11} that
\begin{eqnarray}\label{4.5A}
\pounds_{X}d\eta = d\pounds_{X}\eta = (df)\wedge \eta + f(d\eta).
\end{eqnarray}
Now, making use of (\ref{4.5A}) in \eqref{4.5AB}, we obtain
\begin{eqnarray}\label{4.5X}
&2k(\pounds_{X}\varphi)Z = 4Q\varphi Z + 2(fk - 2\lambda)\varphi Z + k(\eta(Z)Df - (Zf)\xi),
\end{eqnarray}
for any vector field $Z$ on $M$. Next, replacing $\xi$ instead of $Z$ in the last equation and using $\varphi\xi=0$ we have
\begin{eqnarray}\label{4.6}
&2(\pounds_{X}\varphi)\xi =  Df - (\xi f)\xi,
\end{eqnarray}
where we use $k$ is positive. Further, tracing (\ref{1.1}) gives
\begin{eqnarray}\label{4.6B}
k div X = (2n + 1)\lambda - r.
\end{eqnarray}
Let $\omega $ be the volume form of $M$, i.e., $\omega = \eta \wedge (d\eta)^n\neq0$. Taking Lie-derivative of this along the vector field $X$ and using the formula $\pounds_{X}\omega=(div X)\omega$ and equation \eqref{4.5A} yields $div X = (n + 1)f$. By virtue of this, \eqref{4.6B} provides
\begin{eqnarray}\label{4.1}
r = (2n + 1)\lambda - (n + 1)kf.
\end{eqnarray}
Also, Lie-differentiation of $g(\xi,\xi) = 1$ along an arbitrary vector field $X$ on $M$ and by the use of the equations (\ref{1.1}), (\ref{2.7}) yields
\begin{eqnarray}\label{4.2}
k g(\pounds_{X}\xi,\xi) = 2n-\lambda.
\end{eqnarray}
Now, taking Lie-derivative of (\ref{2.1B}) on $X$ and using (\ref{1.1}), (\ref{2.7}) and (\ref{2.11}) we obtain
$k\pounds_{X}\xi = (kf - 2\lambda + 4n)\xi.$
Making use of this in (\ref{4.2}) yields $kf = \lambda - 2n.$ and therefore we have $k\pounds_{X}\xi = (2n-\lambda)\xi$. Next, taking Lie-derivative of $\varphi\xi=0$ along $X$ and using the foregoing equation we get $(\pounds_{X}\varphi)\xi = 0$. In view of this, the Eq$.$ (\ref{4.6}) becomes $Df = (\xi f)\xi$, i.e. $df = (\xi f)\eta$. Exterior derivative of the preceding equation gives $d^2f = d(\xi f)\wedge \eta + (\xi f)d\eta$. Using $d^2 = 0$ in the last equation and then taking the wedge product with $\eta$ we get $(\xi f)\eta \wedge d\eta = 0$. By the definition of contact structure we know that $\eta \wedge d\eta$ is non-vanishing everywhere on $M$. Hence the previous equation provides $\xi f =0$. This implies that $df = 0$, and therefore  $f$ is constant on $M$.
Integrating both sides of  $div X = (n + 1)f$ over $M$ and applying the divergence theorem we get $f = 0$. Since $kf = \lambda - 2n$, it follows that $\lambda = 2n$. Consequently, equation (\ref{4.1}) gives $r = 2n(2n + 1)$. This shows that the scalar curvature is constant.  As $M$ is compact, we may invoke Theorem [WGX] to conclude that $M$ is isometric to $S^{2n + 1}(c)$, where $c = \sqrt[]{\frac{2n(2n + 1)}{r}}$ is the radius of the sphere. Since $r = 2n(2n + 1)$, we have $c = 1$, and hence $M$ is isometric to a unit sphere $S^{2n + 1}$. This completes the proof.~~~~~~~~~~~~~~~~~~~~~~~~~~~~~~~~~~~~~~~~~~~~~~~~~~~~~~~~~~~~~~~~~~~~~~~~~~~~~~~~~~~~~~~~~~$\square$\\

Waiving the compactness assumption and imposing a commutativity condition we have
\begin{theorem}
Let $M^{2n+1}(\varphi,\xi,\eta,g)$, $n>1$, be a $K$-contact manifold with $Q\varphi = \varphi Q$. If $g$ is a $k$-almost Ricci soliton such that $X$ is a contact vector field, then it is trivial and the soliton vector field is Killing.
\end{theorem}
\textbf{Proof: }
As $f$ is constant and $fk = \lambda - 2n$, the equation (\ref{4.5X}) becomes
\begin{eqnarray}\label{4.7}
 k(\pounds_{X}\varphi)Z = 2Q\varphi Z - (\lambda + 2n)\varphi Z,
\end{eqnarray}
for all vector field $Z$ on $M$. Also, from \eqref{2.11} we have $k \pounds_{X}\eta = (\lambda-2n)\eta$. Now, taking Lie-derivative of $\varphi ^2 Z = -Z + \eta(Z)\xi$ along $X$ and then multiplying $k$ on both sides and using the forgoing Eq$.$, we obtain $k \varphi(\pounds_{X}\varphi)Z + k(\pounds_{X}\varphi)\varphi Z = 0$ for any vector field $Z$ on $M$. In view of (\ref{4.7}), the last Eq$.$ becomes
\begin{eqnarray*}
&\varphi Q\varphi Z + Q\varphi^2 Z = (\lambda + 2n)\varphi^2 Z,
\end{eqnarray*}
for any vector field $Z$ on $M$. Since $Q\varphi = \varphi Q$, the last equation reduces to
\begin{eqnarray}\label{4.9}
&QZ = (\frac{\lambda + 2n}{2})Z + (\frac{ 2n - \lambda}{2})\eta(Z)\xi,
\end{eqnarray}
for any vector field $Z$ on $M$. This shows that $(M,g)$ is $\eta$-Einstein. Now, differentiating (\ref{4.9}) along an arbitrary vector field $Y$ on $M$ and using (\ref{2.6}), we get
\begin{eqnarray*}
&(\nabla_YQ)Z = (\frac{Y\lambda}{2})Z - (\frac{Y\lambda}{2})\eta(Z)\xi - (\frac{2n-\lambda}{2})\{g(Z, \varphi Y)\xi + \eta(Z)\varphi Y\},
\end{eqnarray*}
for any vector field $Z$ on $M$. Tracing the foregoing equation over $Y$ and $Z$, respectively, we have $Zr = Z\lambda - (\xi \lambda )\eta(Z)$ and $Zr = n Z(\lambda)$ for any vector field $Z$ on $M$. Since $\xi$ is Killing, $\xi r = 0$. Hence, the last equation provides $\xi \lambda = 0$. Consequently, we have $Zr = Z\lambda$ and $Zr = n Z(\lambda)$ for any vector field $Z$ on $M$. As $n>1$, these two equations imply that $\lambda$ and $r$ are constant. By virtue of (\ref{4.9}), Eq$.$ (\ref{4.7}) reduces to $(\pounds_{X}\varphi)Z = 0$ for any vector field $Z$ on $M$. At this point, we recall Lemma $1$ (cf. \cite{GS}) ``if a vector field $X$ leaves the structure tensor $\varphi$ of the contact metric manifold $M$ invariant, then there exists a constant $c$ such that $\pounds_{X}g = c(g + \eta \otimes \eta)$ " to conclude that
\begin{eqnarray*}
&(\pounds_{X}g)(Y,Z) = c\{g(Y, Z) + \eta(Y)\eta(Z)\},
\end{eqnarray*}
for all vector fields $Y$, $Z$ on $M$. On the other hand, making use of (\ref{4.9}) in the Eq$.$ (\ref{1.1}), we find
\begin{eqnarray*}
&\frac{k}{2}(\pounds_{X}g)(Y,Z) = \lambda g(Y, Z) - S(Y, Z),
\end{eqnarray*}
for all vector fields $Y$, $Z$ on $M$. Comparing the last two equations, we deduce
\begin{eqnarray}\label{4.10}
&\frac{ck}{2}\{g(Y, Z) + \eta(Y)\eta(Z)\} = \lambda g(Y, Z) - S(Y, Z).
\end{eqnarray}
Next, putting $Y = Z = \xi$ in (\ref{4.10}) and using (\ref{2.7}), we get
$ck = 2(\lambda - 2n)$. Further, tracing (\ref{4.10}) yields $ ck(n+1) = (2n+1)\lambda - r$. These two equations together provides $ r = 4n(n+1) - \lambda$. Moreover, using $kf = \lambda - 2n$ in (\ref{4.1}) we have $r = (2n + 1)\lambda - (n + 1)(\lambda - 2n)$. Comparing the last two equations we see that $\lambda = 2n$. Utilizing this in (\ref{4.9}) provides $QY = 2nY$ for any vector fields $Y$ on $M$, i.e. the soliton is trivial. This completes the proof.~~~~~~~~~~~~~~~~~~~~~~~~~~~~~~~~~~~~~~~~~~~~~~~~~~~~~~~~~~~~~~~~~~~~~~~~~~~~~~~~~~~~~~~~~~$\square$\\

For a Sasakian manifold the commutativity condition $Q\varphi = \varphi Q$ holds trivially (e.g., see Blair \cite{Blair}). Thus, we have the following
\begin{corollary}
Let $M^{2n+1}(\varphi,\xi,\eta,g)$, $n>1$, be a Sasakian manifold with $X$ is a contact vector field. If $g$ is a $k$-almost Ricci soliton, then it is trivial and the soliton vector field is Killing.
\end{corollary}

\section{$k$-almost Ricci soliton where  $X = \rho\xi$}
In this section, we shall discuss about some special type of $k$-almost Ricci soliton where the potential vector field $X$ is point wise collinear with the Reeb vector field $\xi$ of the contact metric manifold.

\begin{theorem}
Let  $M^{2n + 1}(\varphi,\xi,\eta,g)$ be a compact $H$-contact manifold. If $g$ represents a non-trivial $k$-almost Ricci soliton with non-zero potential vector field $X$ collinear with the Reeb vector field $\xi$, then $M$ is Einstein and Sasakian.
\end{theorem}
\textbf{Proof:} Since the potential vector field $X$ on $M$ is collinear with the Reeb vector field $\xi$, we have $X = \rho\xi$, where $\rho$ is a non-zero smooth function on $M$ (as $X$ is non zero). Taking covariant derivative along an arbitrary vector field $Y$ on $M$ and using (\ref{2.3}) and \eqref{2.4} we get
\begin{equation}\label{5.0}
\nabla_{Y}X = (Y\rho)\xi - \rho(\varphi Y +\varphi hY).
\end{equation}
By virtue of this the soliton equation (\ref{1.1}) becomes
\begin{eqnarray}\label{5.1}
k(Y\rho)\eta(Z) + k(Z\rho)\eta(Y) - 2k\rho g(\varphi hY,Z) + 2S(Y,Z) = 2\lambda g(Y,Z),
\end{eqnarray}
for all vector fields $Y$, $Z$ on $M$. Replacing $\xi$ instead of $Z$ in (\ref{5.1}) gives
\begin{eqnarray}\label{5.2}
kD\rho +k(\xi\rho)\xi +2(Q\xi - \lambda\xi) =0.
\end{eqnarray}
At this point, putting $Y = Z = \xi$ in (\ref{5.1}) and making use of (\ref{2.4}) yields
\begin{eqnarray}\label{5.3}
 k(\xi\rho) + Trl = \lambda.
\end{eqnarray}
Since $M$ is $H$-contact, the Reeb vector field $\xi$ is an eigenvector of the Ricci operator at each point of $M$, i.e. $Q\xi = (Trl)\xi$. Substituting this in (\ref{5.1}) and then using (\ref{5.3}), we have $kD\rho = k(\xi\rho)\xi$. Since the $k$-almost Ricci soliton is non trivial and $k$ is a positive function, we have $D\rho = (\xi\rho)\xi$. Next, taking covariant derivative along an arbitrary vector field $Y$ on $M$ and using (\ref{2.3}) yields $\nabla_{Y}D\rho = Y(\xi\rho)\xi - (\xi\rho)(\varphi Y +\varphi hY).$
In view of $g(\nabla_{Y}D\rho,Z) = g(\nabla_{Z}D\rho,Y)$, the foregoing equation yields
\begin{eqnarray}\label{5.4}
&Y(\xi\rho)\eta(Z) - Z(\xi\rho)\eta(Y) + 2(\xi\rho)d\eta(Y,Z) = 0,
\end{eqnarray}
for all vector fields $Y$, $Z$ on $M$. Choosing $X$, $Y$ orthogonal to $\xi$ and noting that $d\eta \neq 0$, the last equation provides $\xi\rho = 0$. Hence $\rho$ is constant. Thus, the equation (\ref{5.1}) reduces to
\begin{eqnarray}\label{5.5}
&QZ + (k\rho)h\varphi  Z = \lambda Z,
\end{eqnarray}
for any vector field $Z$ on $M$. Taking the trace of (\ref{5.5}) we obtain $r = (2n+1)\lambda$.
Further, covariant derivative of (\ref{5.5}) along an arbitrary vector field $Y$ on $M$ gives
$$(\nabla_{Y}Q)Z + (k\rho)(\nabla_{Y}h\varphi)Z + \rho (Yk)h\varphi Z = (Y\lambda)Z.$$
Contracting this over $Y$ yields
\begin{eqnarray}\label{5.6}
&\frac{1}{2}(Zr) + \rho ((h\varphi Z)k) + (k\rho)div(h\varphi)Z = (Z\lambda),
\end{eqnarray}
for any vector field $Z$ on $M$. On a contact metric manifold it is known that $div(h\varphi)Z = g(Q\xi,Z) - 2n\eta(Z)$ for all vector field $Z$ on $M$ (see \cite{Blair}). Using $Q\xi = (Trl)\xi$ in the previous equation we have $div(h\varphi)Z = (Trl - 2n)\eta(Z)= |h|^2\eta(Z)$. Hence, equation (\ref{5.6}) reduces to
\begin{eqnarray}\label{5.7}
&\frac{1}{2}(Z r) + \rho ((h\varphi Z)k) + (k\rho)(Trl - 2n)\eta (Z) = (Z\lambda),
\end{eqnarray}
for any vector field $Z$ on $M$. Setting $Z=\xi$ and making use of $r = (2n+1)\lambda$ and (\ref{2.4}) equation (\ref{5.7}) reduces to
\begin{eqnarray}\label{5.8}
\frac{2n-1}{2}(\xi r) - (k\rho)|h|^2 = 0.
\end{eqnarray}
Taking into account (\ref{5.0}) and $X = \rho \xi$ we see that $div (rX) = \rho(\xi r) +r(\xi \rho) = \rho(\xi r)$, where we have also used $tr(h\varphi)=0$. Using this equation in  (\ref{5.8}) gives $ \frac{2n-1}{2}div (rX) = (k\rho^2)|h|^2$. Integrating this over $M$ and using divergence theorem we obtain
$$ \int {k\rho^2  \left| h \right|^2} dM = 0.$$
Since the soliton is non-trivial with non-zero potential vector field $X$ and $k$ being positive, the foregoing equation implies $h = 0$ and hence $M$ is $K$-contact. Therefore, equation (\ref{5.5}) shows that $QZ = \lambda Z$. Using (\ref{2.7}) it follows that $\lambda = 2n$. Thus, $M$ is Einstein with Einstein constant $2n$. So, we can apply the result of Boyer-Galicki \cite{BG} which states that ``any compact $K$-contact Einstein manifold is Sasakian" to conclude the proof.~~~~~~~~~~~~~~~~~~~~~~~~~~~~~~~~~~~~~~~~~~~~~~~~~~~~~~~~~~~~~~~~~~~~~~~~~~~~~~~~~~~~~~~~~~$\square$\\

For a $K$-contact manifold it is known that $Q\xi = 2n\xi$. Thus, from the above theorem we have the following
\begin{corollary}
If a complete $K$-contact metric represents a non-trivial $k$-almost Ricci soliton with non-zero potential vector field $X$ collinear with the Reeb vector field $\xi$, then it is Einstein and Sasakian.
\end{corollary}

Next, replacing the ``compact $H$-contact" of the previous theorem by the commutativity condition $Q\varphi = \varphi Q$ we prove
\begin{theorem}
Let $M^{2n + 1}(\varphi,\xi,\eta,g)$ be a contact metric manifold satisfying $Q\varphi = \varphi Q$. If $g$ represents a non-trivial $k$-almost Ricci soliton with nonzero potential vector field $X$ collinear with the Reeb vector field $\xi$, then $M$ is Einstein and $K$-contact. In addition, if $M$ is complete, then it is compact Sasakian.
\end{theorem}
\textbf{Proof :} The commutativity condition $Q\varphi = \varphi Q$ together with (\ref{2.4}) and $\varphi\xi=0$ shows that $Q\xi = (Tr l)\xi$. Further, since the potential vector field $X$ is collinear with the Reeb vector field $\xi$, from equations (\ref{5.0}) to (\ref{5.8}) are also valid here. Now we replace $Z$ by $\varphi Z$ in (\ref{5.5}) and using $h\xi=0$ we get
\begin{eqnarray}\label{5.9}
	&Q\varphi Z - (k\rho)hZ = \lambda \varphi Z.
\end{eqnarray}
On the other hand, operating (\ref{5.5}) by $\varphi$ and using $h\varphi=-\varphi h$ we obtain
\begin{eqnarray}\label{5.10}
	&\varphi QZ + (k\rho)hZ = \lambda \varphi Z.
\end{eqnarray}
Adding (\ref{5.9}) and (\ref{5.10}) along with $Q\varphi = \varphi Q$ gives $Q\varphi Z = \lambda \varphi Z$. Therefore, replacing $Z$ by $\varphi Z$ in the last equation and using $Q\xi = (Tr l)\xi$, we deduce $QZ = \lambda Z + (Tr l - \lambda )\eta(Z)\xi $. Since $\lambda = Tr l$ (follows from \eqref{5.3}, as $\xi \rho=0$), the foregoing equation implies $QZ = \lambda Z$ and hence $M$ is Einstein. Consequently, the scalar curvature $r$ and $\lambda$ are constant. Thus, from (\ref{5.8}), it follows that $(k\rho)|h|^2 = 0.$ Since $k$ is positive and the soliton vector field $X$ is non-zero, we can conclude that $h = 0$, and hence $M$ is $K$-contact. From these, we see that $M$ is $K$-contact and Einstein with Einstein constant $2n$. Now, if $M$ is complete then applying Myers' theorem $M$ becomes compact. Finally, using Boyer-Galicki's theorem \cite{BG} we conclude the proof.~~~~~~~~~~~~~~~~~~~~~~~~~~~~~~~~~~~~~~~~~~~~~~~~~~~~~~~~~~~~~~~~~~~~~~~~~~~~~~~~~~~~~~~~~~$\square$

\noindent\\
\textbf{Acknowledgments: }The authors are very much thankful to the reviewer for his/her comments. The author D. S. Patra is financially supported by the Council of Scientific and Industrial Research, India (grant no. 17-06/2012(i)EU-V).

\noindent
$^1$
Department of Mathematics, \\
Chandernagore College\\
Hooghly: 712 136 (W.B.), INDIA\\
E-mail: aghosh\_70@yahoo.com\\

\noindent
$^2$
Department of Mathematics, \\
Jadavpur University,  \\
188. Raja S. C. Mullick Road,\\
Kolkata:700 032, INDIA \\
E-mail: dhritimath@gmail.com

\end{document}